\documentclass{amsart}
\usepackage{amsmath,amssymb}
\usepackage[all]{xy}

\theoremstyle{plain}
\newtheorem{introtheorem}{Theorem}
\newtheorem{theorem}{Theorem}[section]
\newtheorem{proposition}[theorem]{Proposition}
\newtheorem{lemma}[theorem]{Lemma}
\newtheorem{corollary}[theorem]{Corollary}

\theoremstyle{definition}

\newtheorem{example}[theorem]{Example}

\theoremstyle{remark}
\newtheorem{remark}[theorem]{Remark}

\newcommand{\secref}[1]{Section~\ref{#1}}
\newcommand{\thmref}[1]{Theorem~\ref{#1}}
\newcommand{\propref}[1]{Proposition~\ref{#1}}
\newcommand{\lemref}[1]{Lemma~\ref{#1}}
\newcommand{\corref}[1]{Corollary~\ref{#1}}

\newcommand{\exref}[1]{Example~\ref{#1}}

\def\Z{{\mathbb Z}}
\def\Q{{\mathbb Q}}
\def\C{{\mathbb C}}

\def\map{\mathrm{map}}
\def\Der{\mathrm{Der}}
\def\Hom{\mathrm{Hom}}
\def\Rel{\mathrm{Rel}}
\def\cat0{\mathrm{cat}_0}
\def\rank{\mathrm{rank}}
\def\dim{\mathrm{dim}}

\def\im{\mathrm{im}}

\begin{document}

\title[Rank of the Fundamental Group of   Function Spaces]
{Rank of the Fundamental Group of  a Component of a Function Space}

\author{Gregory  Lupton}

\address{Department of Mathematics,
          Cleveland State University,
          Cleveland OH 44115}

\email{G.Lupton@csuohio.edu}

\author{Samuel Bruce Smith}

\address{Department of Mathematics,
  Saint Joseph's University,
  Philadelphia, PA 19131}

\email{smith@sju.edu}

\date{\today}



\begin{abstract}  We compute the rank of the fundamental group of 
an arbitrary connected component  of   the space $\map(X, Y)$  
 for $X$ and $Y$ nilpotent 
 CW complexes with $X$ finite.  For the general component 
 corresponding to a homotopy class  $f \colon X \to 
 Y,$ we give a  formula directly computable  from the Sullivan model for 
 $f$.  For the component of the constant map,  
  our  formula expresses the rank in terms of classical invariants of 
  $X$ and $Y$.   
  Among other applications and calculations, we obtain the following: 
  Let  $G$ be  a compact simple Lie group with maximal torus $T^n$.
  Then $\pi_1(\map(S^2, G/T^n; f))$ is a finite group if and only if 
  $f \colon S^2 \to G/T^n$ is essential. 
\end{abstract}

\maketitle

\section{Introduction}\label{sec:intro}
Given a finitely generated abelian group $G$, the {\em rank} of $G$ is 
defined to be the cardinality of a basis  of the free abelian group 
$G/T$ where $T$ is the torsion subgroup of $G$.  This notion  
can be extended to finitely generated nilpotent groups 
$\Gamma$ by considering the 
central series 
$ \Gamma = \Gamma_0 \supset \Gamma_1 \supset \cdots \supset 
\Gamma_n = \{1 \}$
where $\Gamma_n = [\Gamma_{n-1}, \Gamma]$.  Specifically, the {\em rank}
of $\Gamma$ is defined to be the sum of the ranks of the finitely 
generated abelian quotients $\Gamma_{k-1}/\Gamma_{k}$ for $k = 1, 
\ldots, n$. Alternately,    
$\rank(\Gamma) = \dim_\Q(\Gamma_\Q)$ where
$\Gamma_\Q$ is the rationalization of the   nilpotent group $\Gamma$
as described   in \cite[Ch. 1]{H-M-R}. 

Given  connected spaces  $X$ and $Y$ let $\map(X, Y)$ denote the 
 function space of  all (not-necessarily
based) maps from $X$ to $Y$ with the compact-open topology.  Given a 
particular map $f \colon X \to Y$, write
$\map(X,Y;f)$ for the connected component of $\map(X, Y)$ containing 
$f$, that is, the space of maps    
freely homotopic to $f$. Under reasonable hypotheses on $X$ and $Y$, 
the fundamental groups of the connected components of $\map(X, 
Y)$ are  finitely generated nilpotent groups.  
To be precise, 
  say a space $X$ is  {\em 
nilpotent} if $\pi_1(X)$ is a nilpotent group and the action of 
$\pi_1(X)$ on the higher homotopy groups of $X$ is a nilpotent 
action.   If   $X$ is a finite CW complex and $Y$ is a  nilpotent CW complex
of finite type
then the components of  $\map(X, Y)$ are themselves  nilpotent CW 
complexes of finite type
(\cite{Mil} and \cite[Th.II.2.5]{H-M-R}).  
Thus with these hypotheses 
$\pi_1\big{(}\map(X, Y;f)\big{)}$ is a finitely generated 
nilpotent group. 
Our purpose in this paper is to  give a formula for the 
rank of this  group  expressed in terms of accessible invariants of the 
map $f \colon X \to Y.$

In  \cite{L-S}, we 
describe the higher rational homotopy groups of $\map(X,Y;f).$ Our
description uses Sullivan minimal models, and is in terms of the homology
of chain complexes of derivations. Specifically, for each $n\geq
2$ we construct an isomorphism 
$$\pi_n\big(\map(X,Y;f)\big)\otimes\Q \cong
H_n\big(\Der(\mathcal{M}_Y, \mathcal{M}_X; \mathcal{M}_f)\big)$$
where $\mathcal{M}_f \colon \mathcal{M}_Y \to \mathcal{M}_X$ is the 
Sullivan model of the map $f \colon X \to Y$ and
\linebreak
$\Der_*(\mathcal{M}_Y, \mathcal{M}_X; \mathcal{M}_f)$ is a differential
graded vector space of (generalized) algebra derivations (see  
\secref{sec:derivations} below for precise definitions).
In \cite{L-S}---and its companion \cite{L-S2}, in which we give
corresponding results using differential graded Lie algebra minimal models---our
main goal was to establish a framework within which we could study
evaluation subgroups of $\pi_n(Y)\otimes\Q$, for $n \geq 2$.
Now within that framework, we may also compute
$H_1\big(\Der(\mathcal{M}_Y, \mathcal{M}_X; \mathcal{M}_f)\big)$,
which is a rational homotopy invariant of the map $f$. Heuristically, one
would expect this vector space to be related to rationalization of the fundamental group
$\pi_1\big(\map(X,Y;f)\big)$. However,  the latter group is generally non-abelian,
 and so the above isomorphism obviously cannot be extended in the
 na\"{\i}ve way.  
As the main result of this paper, we establish   
\begin{introtheorem}\label{thm:general component}
Let $X$ and $Y$ be connected nilpotent CW complexes of finite type with 
$X$ finite. Let $f \colon X \to Y$ be a given map. 
Then 
$$\rank\big(\pi_1\big(\map(X,Y;f)\big)\big) = 
\dim_\Q\big(H_1\big(\Der(\mathcal{M}_Y, \mathcal{M}_X; \mathcal{M}_f)\big)\big).$$
\end{introtheorem}

In the special case of a null component, that is,  $\map(X,Y;0)$, we obtain
an identification of the rank 
of the fundamental group directly in 
terms of classical invariants of $X$ and $Y$. 
Let $b_n(X)$ denotes the $n$th
Betti number of $X$.  If $n\geq2$ let $\rho_n(Y)$ denote the
integer $\dim_\Q(\pi_n(Y)\otimes\Q)$---sometimes called the $n$th
\emph{Hurewicz number} of $Y$. We have 
\begin{introtheorem}\label{thm:null component}
Let $X$ and $Y$ be nilpotent spaces, with $X$ a finite complex of
dimension $N$.  Then  
$$\rank\big(\pi_1\big(\map(X,Y;0)\big)\big) = 
\rank\big(\pi_1(Y)\big) +  \sum_{n=2}^{N+1}
\rho_n(Y)\cdot b_{n-1}(X)  .$$
\end{introtheorem}
\begin{remark} 
So far as the rational homotopy groups of function spaces are 
concerned, these results form an essentially 
complete complement to the results of \cite{L-S}.  Of course, 
they leave open the more general problem of  
describing   the 
full structure of the rationalized fundamental group 
 $\pi_1\big(\map(X,Y;f)\big)_\Q$. 
In principle, 
 such a  description  is already available 
 via   the   Sullivan models 
 for $\map(X,Y;f)$ described by Haefliger \cite{Hae} and 
 Brown-Sczcarba \cite{B-Sz}. 
 In practice,  
 however, there are  so    many technical issues involved in first
 assembling these   models, and then extracting from them the necessary
 information, that these identifications are not  helpful for our 
 purposes.  By
 focussing purely on the rank of $\pi_1\big(\map(X,Y;f)\big)$, we
 are able to give a ``closed form" description that proceeds
 directly from the ordinary Sullivan model of the map $f \colon X
 \to Y$ and is relatively easy to handle in practice. 
 We emphasize that our results and methods here are entirely 
 independent of those of \cite{Hae} and \cite{B-Sz}
 At the end of 
 the paper,  we consider some particular cases
 where further structure  of the rationalized 
 fundamental group can be described. 
 \end{remark}
 \begin{remark}
 \thmref{thm:null component}  extends to the fundamental group an isomorphism for higher 
homotopy groups of the form 
 $$\pi_n\big( \map(X, Y;0) \big) \otimes \Q 
 \cong  \bigoplus_{k \geq n} H^{k-n}(X, \pi_{n}(Y) \otimes \Q)$$ 
 $n \geq 2$.   The isomorphism for $n \geq 
 2$  has appeared in a number of places 
 (see \cite{Sm94, B-Z97, L-S}).  
It is worth noting that the assertions of both \cite[Lem.3.2]{Sm94} 
and  \cite{B-Z97}
 need to be adjusted so as to exclude the
$\pi_1$ case, even under the hypothesis that $Y$ is $1$-connected.
\end{remark}

We give several interesting computations using our results.  In 
\exref{ex:free loop}, we compute the rank of the fundamental group
of any arbitrary component  of the free loop space for  $Y$. 
In \thmref{thm:F_0-space}, we specialize our  results  to  a 
  class of 
 formal spaces which includes all homogeneous spaces $G/H$ of equal 
 rank Lie pairs.  We express the 
 rank  of an arbitrary component corresponding to $f$ in this case 
 directly in terms of 
 the map induced by $f$ on rational cohomology.
 As a particular consequence, in \exref{ex:G/T} we show the following:
 Let $G$ be a compact simple Lie group of rank $n \geq 2,$ $T^n 
 \subseteq G$ a maximal torus and $f 
 \colon S^2 \to G/T^n$ any map.  Then $\pi_1\big(\map(S^2, 
 G/T^n;f)\big)$ is finite -- that is, has rank $0$ -- if and only if 
 $f$ is essential.

We now give an outline of the paper.  \secref{sec:example}
consists of a concrete example that illustrates a 
component of a function space  may have fundamental group 
that is non-abelian after rationalization.  The
example suggests that rationalized fundamental groups of function
space components form a rich subject for investigation.
 In \secref{sec:minimal models}
we review some properties of minimal models of nilpotent spaces and  in 
\secref{sec:derivations} we introduce the framework of derivation chain complexes  
in preparation for the  proof of   the  main results in 
\secref{sec:proofs}.  In  
\secref{sec:consequences}   we deduce several
consequences of  our results and give some further examples. 

We end this section by fixing some notation and terminology for the sequel. 
We denote
the trivial group by $\{1\}$, the trivial vector space by $\{0\}$,
and the space consisting of a single point by $\{*\}$.  We also
use $0$ to denote a trivial homomorphism of groups or vector
spaces, or a trivial map of spaces. 
We use $H(g)$, respectively $g_\#$, to denote a homomorphism induced
by $g$ on either cohomology or homology, respectively homotopy.
Whether we intend cohomology or homology, or the use of particular
coefficients, will be clear from context.

We assume basic familiarity with the localization of 
nilpotent groups and spaces, as discussed in
\cite{H-M-R}. In particular, we recall that a nilpotent space $Y$ 
(respectively, group $G$)  
admits a {\em rationalization} which we denote $e_Y \colon Y \to Y_\Q$ 
(respectively, $e_G \colon G \to G_\Q$).  
Recall that if $G$ is abelian then $G_\Q \cong G 
\otimes \Q.$ 
If $f \colon X \to Y$ is a map into a
nilpotent space, we write  $f_\Q \colon X \to Y_\Q$ for the 
compositon $  e_Y \circ f$ and likewise for a group homomorphism. 
A map of spaces $g \colon Y \to Z$ induces a 
map  $g_* \colon \map(X, Y; f)
\to \map(X, Z; g\circ f)$ induced by post-composition with $g$. 
By  \cite[Th.II.3.1]{H-M-R}, 
if $f \colon X \to Y$ is a map between spaces $X$ and $Y$ satisfying  the 
hypotheses in the theorems above,  
 then the induced map
$(e_Y)_* \colon \map(X, Y;f) \to \map(X, Y_\Q; f_\Q)$ 
is a rationalization of the nilpotent space $\map(X, Y;f)$.   In particular, 
under these hypotheses  
\begin{equation}\label{eq:map rationalized} 
    \pi_n\big( \map(X, Y;f) \big)_\Q \cong \pi_n\big( \map(X, Y_\Q; 
f_\Q) \big)\end{equation}
for all $n \geq 1.$

\section{An Example}\label{sec:example}

We begin with a concrete example.  It illustrates that
$\map(X,Y;0)$ may have fundamental group that is non-abelian after
rationalization, even though both $X$ and $Y$ are simply
connected. In the course of our discussion, we establish a basic
ingredient of our further developments.  This section does not
require any familiarity with minimal models.

The following is a particular phrasing of a well-known result.

\begin{lemma}\label{lem:pullback}
Suppose given a pullback square of topological spaces and maps
$$\xymatrix{A \ar[r]^{\bar k} \ar[d]_{\bar h} & B \ar[d]^{h}\\
C \ar[r]_{k} & D.}$$
If $f \colon X \to A$ is any map, then
$$\xymatrix{\map(X, A; f) \ar[r]^-{({\bar k})_*} \ar[d]_-{({\bar h})_*} & \map(X, B; {\bar k}\circ f) \ar[d]^-{h_*}\\
\map(X, C; {\bar h}\circ f) \ar[r]_-{k_*} & \map(X, D; k\circ{\bar
h}\circ f = h\circ{\bar k}\circ f).}$$
is a pullback square.
\end{lemma}

\begin{proof}
Suppose given based maps $\beta\colon Z \to \map(X, B; {\bar k}
\circ f)$ and $\gamma\colon Z \to \map(X, C;$ ${\bar h} \circ f)$
such that $h_*\circ\beta = k_*\circ\gamma$.  We must show that
there exists a unique based map $\alpha\colon Z \to \map(X, A; f)$
that satisfies $({\bar k})_*\circ\alpha = \beta$ and $({\bar
h})_*\circ\alpha = \gamma$.  To this end, let $b\colon Z\times X
\to B$ and $c\colon Z\times X \to C$ be the adjoints of $\beta$
and $\gamma$, respectively.  Because of the choice of components
for each function space, the existence and uniqueness of $\alpha$
is equivalent to the existence and uniqueness, respectively, of
its adjoint $a$ in the following commutative diagram:
\begin{displaymath}
\xymatrix{ Z\times X \ar@/_/[ddr]_{c} \ar@{.>}[dr]_(0.6){a}
\ar@/^/[drr]^{b}&&\\ &A\ar[r]^{\bar k} \ar[d]^{\bar h}&B
\ar[d]^{h}\\ &C \ar[r]_{k}&D\\}
\end{displaymath}
The result follows.
\end{proof}

We say a fibration $p \colon E \to B$ is {\em principal} if $p$
is obtained as a pullback
$$\xymatrix{E \ar[r]^{\bar k} \ar[d]_{p} & PK \ar[d]^{q}\\
B \ar[r]_{k} & K,}$$
where $q\colon PK \to K$ is the usual path fibration over some space 
$K$ and $k\colon
B \to K$ is some map.

\begin{corollary}\label{cor:principal}
Suppose given a map $f \colon X \to E$ and a principal fibration
$p \colon E \to B$.  Then the induced map $p_*\colon \map(X,E;f)
\to \map(X, B; p\circ f)$ is a principal fibration.
\end{corollary}

\begin{proof}
 Applying \lemref{lem:pullback} to the diagram above, we obtain a
pullback
$$\xymatrix{\map(X, E; f) \ar[r]^-{({\bar k})_*} \ar[d]_-{p_*} & \map(X, PK; {\bar k}\circ f) \ar[d]^-{q_*}\\
\map(X, B; p\circ f) \ar[r]_-{k_*} & \map(X, K; k\circ p\circ
f).}$$
Since $PK$ is contractible, $\map(X, PK)$ is connected and 
we have a natural identification $P\map(X, K; k\circ p\circ
f) = \map(X, PK).$
\end{proof}

Now take $X = \C P^2$ and $Y$ the total space of the principal
fibration $K(\Z,5) \to Y \to K(\Z\times\Z,3)$ whose $k$-invariant
is $k\colon K(\Z\times\Z,3) \to K(\Z,6)$, corresponding to the
cup-product $x\cup y \in H^6(K(\Z\times\Z,3);\Z)$.  Here, $x$ and
$y$ denote the generators in $H^3(K(\Z\times\Z,3);\Z) \cong \Z
\times \Z$.  Apply \corref{cor:principal} to the principal
fibration $K(\Z,5) \to Y \to K(\Z\times\Z,3)$, and the trivial map
$0\colon X \to Y$.  We obtain a principal fibration $p_*\colon
\map(X,Y;0) \to \map(X, K(\Z\times\Z,3); 0)$.  Since $\map(X,
PK(\Z\times\Z,3); 0)$ is contractible, we have a homotopy pullback
$$\xymatrix{\map(X, Y; 0) \ar[r]^-{0} \ar[d]_-{p_*} & {\{*\}} \ar[d]^-{0} \\
\map(X, K(\Z\times\Z,3); 0) \ar[r]_-{k_*} & \map(X, K(\Z,6); 0)}$$
and consequently a fibre sequence
\begin{equation}\label{eq:fibre sequence}
\xymatrix{\map(X, Y; 0) \ar[r]^-{p_*} & \map(X, K(\Z\times\Z,3);
0) \ar[r]^-{k_*} & \map(X, K(\Z,6); 0)  .}
\end{equation}
From \cite[Th.2]{Th} (see also \cite{Hae}), we have homotopy
equivalences
$$\map(X, K(\Z\times\Z,3); 0) \simeq K(\Z,1)\times K(\Z,1)\times
K(\Z,3)\times K(\Z,3),$$
and
$$\map(X, K(\Z,6); 0) \simeq K(\Z,2)\times K(\Z,4)\times
K(\Z,6).$$
In particular, we have $H^2(\map(X, K(\Z\times\Z,3); \Q) \cong
H^2(\map(X, K(\Z,6); 0); \Q) \cong \Q$.  Now a careful reading of
\cite[Sec.1.2]{Hae} shows that
$$H(k_*)\colon H^2(\map(X, K(\Z,6); 0); \Q) \to H^2(\map(X,
K(\Z\times\Z,3); \Q)$$
is an isomorphism. We now show from this that $\pi_1\big(\map(X,
Y; 0)\big)_\Q$ is non-abelian. From the long exact sequence
induced by (\ref{eq:fibre sequence}) in rational homotopy, the
Serre exact sequence induced by (\ref{eq:fibre sequence}) in
rational homology, and the rational Hurewicz homomorphism between
them, we obtain a commutative ladder as follows:
$$
\xymatrix{\{1\} \ar[r]^-{(k_*)_\#} & \Q
\ar[r]^-{\partial_\#}\ar[d]_{h}^-{\cong} & \pi_1\big(\map(X, Y;
0)\big)_\Q
\ar[r]^-{(p_*)_\#}\ar[d]_{h} & \Q\oplus\Q \ar[r] \ar[d]_{h}^-{\cong} & \{1\}\\
\Q \ar[r]^-{H(k_*)}_-{\cong} & \Q \ar[r]^-{0} & H_1(\map(X, Y; 0);
\Q) \ar[r]^-{H_*(p_*)} & \Q\oplus\Q \ar[r] & \{0\} }
$$
From the top row, $\pi_1\big(\map(X, Y; 0)\big)_\Q$ has rank $3$
(as a nilpotent group), but from the bottom row $H_1(\map(X, Y;
0); \Q)$ has rank (or dimension) $2$.  Therefore, the Hurewicz
homomorphism $h \colon \pi_1\big(\map(X, Y; 0)\big)_\Q \to
H_1(\map(X, Y; 0); \Q)$ is not an isomorphism, and hence
$\pi_1\big(\map(X, Y; 0)\big)_\Q$ is not abelian.

\section{Minimal Models in the Non-Simply Connected Setting}%
\label{sec:minimal models}

We assume familiarity with rational homotopy theory using the DG
(differential graded) algebra minimal models introduced by
Sullivan. Our main reference for this material is \cite{F-H-T},
although that book restricts to the simply connected case.
References that treat the non-simply connected case in some detail
include \cite{Bo-Gu, G-M}. Here we briefly review some of this
material in the nilpotent setting, and take this opportunity to
establish some notation.

In general, we use the standard notation and terminology for
minimal models as  in \cite{F-H-T}.  The basic facts that we
rely on are as follows: Each nilpotent space $Y$ has a Sullivan
minimal model $(\mathcal{M}_{Y}, d_Y)$ in the category of
nilpotent DG algebras over $\Q$. This DG algebra is unique up to
isomorphism and is of the form $\mathcal{M}_{Y} = \Lambda W$, a
free graded commutative algebra generated by a positively graded
vector space $W$ of finite type. The differential $d_Y$ is
decomposable, in that $d_Y(W) \subseteq \Lambda^{\geq 2} W$, and
satisfies a certain ``nilpotency" condition. A map $f\colon X \to
Y$ of nilpotent spaces has a Sullivan minimal model which is a DG
algebra map $\mathcal{M}_f \colon \mathcal{M}_Y \to
\mathcal{M}_X$. The Sullivan minimal model is a complete rational
homotopy invariant for a space or a map. Since the minimal model
is determined by the rational homotopy type, the minimal models of
$Y_\Q$ and $Y$, and more generally those of $f_\Q$ and $f$, agree.
The homomorphism of rational homotopy groups induced by a map
$f\colon X \to Y$ of nilpotent spaces may be identified with the
homomorphism induced by $\mathcal{M}_f$ of the (quotient) modules
of indecomposables $Q(\mathcal{M}_f) \colon Q(\mathcal{M}_Y) \to
Q(\mathcal{M}_X)$.

We now recall the structure of the minimal model of a nilpotent space $Y$.   
By \cite[Th.2.9]{H-M-R}, the  Postnikov
decomposition of $Y$  admits a
{\em principal refinement} at each stage. Precisely,   $p_r \colon Y^{(r)} \to
Y^{(r-1)}$, the $r$th stage of the Postnikov decomposition of
$Y$, factors into a finite  sequence of principal fibrations
\begin{equation}\label{eq:principal factors Y1f}
Y^{(r)} = Y^{(r)}_{c_r} \to Y^{(r)}_{c_r - 1} \to \cdots \to
Y^{(r)}_{1}\to Y^{(r)}_{0} = Y^{(r-1)},
\end{equation}
each induced from a path-loop fibration by a $k$-invariant of the
form $k^r_j \colon  Y^{(r)}_{j-1} \to K(G\,^r_j, r+1)$.  

The
minimal model $\mathcal{M}_Y$ of $Y$ may be constructed by a
sequence of so-called \emph{elementary extensions} in a way that
mirrors this principal refinement. 
Thus, for example, consider the \emph{$1$-minimal model} of $Y$,
that is, the sub-DG algebra $\mathcal{M}_{Y^{(1)}}$ of
$\mathcal{M}_Y$ generated in degree $1$.   Let $V\,^1_j$ and
$\overline{V}\,^1_j$ denote vector spaces isomorphic to
$G\,^1_j\otimes\Q$ and concentrated in degree $2$ and $1$
respectively. Then the $1$-minimal model $\mathcal{M}_{Y^{(1)}} =
\Lambda(W^{(1)}, d)$ is $c_1$-stage in the sense that $W^{(1)} =
\oplus_{j=1}^{c_1} \overline{V}\,^1_j$, with
$d(\overline{V}\,^1_1) = 0$ and $d(\overline{V}\,^1_j) \subseteq
\Lambda(\overline{V}\,^1_1\oplus
\cdots\oplus\overline{V}\,^1_{j-1})$ for $j = 2, \dots, c_1$. For
each stage of the $1$-minimal model we have an elementary K-S
extension
\begin{equation}\label{eq:K-S extn}
\mathcal{M}_{Y^{(1)}_{j-1}} \to \mathcal{M}_{Y^{(1)}_{j}} \to
(\Lambda \overline{V}\,^1_j,d=0)
\end{equation}
that is a K-S model of the principal fibration $K(G^1_j, 1) \to
Y^{(1)}_{j} \to Y^{(1)}_{j-1}$.  This extension is elementary in
the sense that $\mathcal{M}_{Y^{(1)}_{j}} =
(\mathcal{M}_{Y^{(1)}_{j-1}}\otimes\Lambda \overline{V}\,^1_j, D)$
and $D(\overline{V}\,^1_j) \subseteq \mathcal{M}_{Y^{(1)}_{j-1}}$.
We write $V\,^1_j = \langle v^1_{j,1}, \dots, v^1_{j,n^1_j}
\rangle$ and $\overline{V}\,^1_j = \langle \overline{v}^1_{j,1},
\dots, \overline{v}^1_{j,n^1_j} \rangle$.  Then the extension
(\ref{eq:K-S extn}) is minimal in the sense that, for each $i$,
$D(\overline{v}\,^1_{j,k}) = \xi^1_{j,k}$ with each $\xi^1_{j,k}$
decomposable in $\mathcal{M}_{Y^{(1)}_{j-1}}$.  In fact, since we
only have generators of degree $1$ so far, each $\xi^1_{j,k}$ is
of homogeneous length $2$. The relation between $k$-invariants of
the principal fibrations of (\ref{eq:principal factors Y1f}) and
the elementary extensions (\ref{eq:K-S extn}) is as follows: Each
$\xi^1_{j,k}$ is a cycle that represents a class in
$H^2(\mathcal{M}_{Y^{(1)}_{j-1}})$. On the other hand, the
$k$-invariant $k^1_j$ has minimal model $\mathcal{M}_{k^1_j}
\colon (\Lambda V\,^1_j,d=0) \to \mathcal{M}_{Y^{(1)}_{j-1}}$,
with $\mathcal{M}_{k^1_j}(v^1_{j,k}) = \xi^1_{j,k}$ for each $j$
and $k$.  Finally, we remark that the $1$-minimal model
$\mathcal{M}_{Y^{(1)}}$ is a complete rational invariant for the
fundamental group of $Y$, in the case in which $Y$ is a nilpotent
group.  In particular, the rank of $\pi_1(Y)$ equals the number of
generators of $\mathcal{M}_{Y^{(1)}}$.  Other aspects of
$\pi_1(Y)$ are determined in different, often less direct, ways by
$\mathcal{M}_{Y^{(1)}}$.  The nilpotency class of $\pi_1(Y)_\Q$,
for example, is determined as the smallest $c_1$ for which
$\mathcal{M}_{Y^{(1)}}$ is $c_1$-stage in the above sense.

A similar situation pertains for the higher dimensional parts of
the minimal model $\mathcal{M}_Y$.  We extend the preceding
notation as follows:  The $r$-minimal model of $Y$, written
$\mathcal{M}_{Y^{(r)}}$, is the sub-DG algebra of $\mathcal{M}_Y$
generated in degrees $\leq r$.  In fact it is a minimal model for
the $r$th stage of the Postnikov decomposition $Y^{(r)}$ of $Y$.
Let $V\,^r_j$ and $\overline{V}\,^r_j$ denote vector spaces
isomorphic to $G\,^r_j\otimes\Q$ and concentrated in degree $r+1$
and $r$ respectively.  Corresponding to (\ref{eq:principal factors
Y1f}) we have $c_r$ minimal, elementary extensions
\begin{equation}\label{eq:K-S extnr}
\mathcal{M}_{Y^{(r)}_{j-1}} \to \mathcal{M}_{Y^{(r)}_{j}} \to
(\Lambda \overline{V}\,^r_j,d=0)
\end{equation}
Write $V\,^r_j = \langle v^r_{j,1}, \dots, v^r_{j,n^r_j} \rangle$
and $\overline{V}\,^r_j = \langle \overline{v}^r_{j,1}, \dots,
\overline{v}^r_{j,n^r_j} \rangle$.  For each $k$,
$D(\overline{v}\,^r_{j,k}) = \xi^r_{j,k} \in
\mathcal{M}_{Y^{(r)}_{j-1}}$.  In the general case, each
$\xi^1_{j,k}$ need not be of homogeneous length but is
decomposable.  Each $\xi^r_{j,k}$ is a cycle that represents a
class in $H^{r+1}(\mathcal{M}_{Y^{(r)}_{j-1}})$. The $k$-invariant
$k^r_j$ has minimal model $\mathcal{M}_{k^r_j} \colon (\Lambda
V\,^r_j,d=0) \to \mathcal{M}_{Y^{(2)}_{j-1}}$, and we have
$\mathcal{M}_{k^r_j}(v^r_{j,k}) = \xi^r_{j,k}$ for each $j$ and
$k$.  Roughly speaking, the way in which generators of degree $1$
are involved with the higher degree generators of
$\mathcal{M}_{Y}$ corresponds to the way in which the fundamental
group of $Y$ is involved with the topology of $Y$.  For instance,
a differential from (\ref{eq:K-S extnr}) of the form
\begin{align*}
D(\overline{v}\,^r_{j,k}) = \sum\ c^k_{s, t, p,q}
\,\overline{v}\,^1_{s,t}\overline{v}\,^r_{p,q} & + \text{quadratic
terms not involving $\mathcal{M}_{Y^{(1)}}$ } \\
&+ \text{length $\geq 3$ terms},
\end{align*}
with at least one coefficient $c^k_{s, t, p,q}$ non-zero, occurs
when $\pi_1(Y)_\Q$ acts non-trivially on $\pi_r(Y)_\Q$.

\section{Derivations and  Function Space Components }
\label{sec:derivations}

In this section, we describe the framework of chain 
complexes of derivations mentioned in the introduction.  We then  construct a 
commutative ladder linking  the long exact homotopy   
sequence of a  principal fibre sequence of function spaces  to 
a  long exact homology sequence within this framework. 
 
Fix  two  DG algebras $(A,d_A)$ and $(B,d_B)$ 
and a  DG algebra
map $\phi \colon A \rightarrow B$ between them.  
Define a \emph{$\phi$-derivation}  to be 
a linear map $\theta
\colon A \rightarrow B$ that \emph{reduces 
degree by $n$} and
satisfies the derivation law $$\theta (xy) = 
\theta(x)\phi(y) +
(-1)^{n|x|} \phi(x) \theta(y).$$  
Let $\Der_{n}(A, B; 
{\phi})$ denote the vector
space of $\phi$-derivations of degree $n \geq 0.$ 
Define a linear map $\delta \colon 
\Der_{n}(A, B; {\phi}) \to
\Der_{n-1}(A, B; {\phi})$ by $\delta(\theta) 
= d_B \circ \theta  -
(-1)^{|\theta|} \theta \circ d_A$.   A 
standard check now shows
that $\delta^2 = 0$ and thus 
$(\Der_{*}(A, B; {\phi}),
\delta)$ is a chain complex.   

The construction is the obvious    generalization of  the standard
DG Lie  algebra of derivations of a DG Lie algebra complex, which we 
would denote  here by
$\Der_{*}(B, B; 1)$.  Of course, in the general case of interest here there 
is no  bracket. In another  special case, when  $(A, d_A) = (\Lambda V, 
0),$
we have 
$ \Der_n(\Lambda V, B; 0), \delta) \cong \Hom_n(V,B)$
where $\Hom_n$ denotes the space of homomorphisms that reduce degree $n$. 
While it is sometimes convenient to truncate these complexes 
taking only cycles in degree $1$,  
we do not take this convention here.  
To reduce notation,  we will   suppress the differential
when writing the homology of a DG algebra. Also, we will only write 
the outermost degree in the homology of a DG space.   In particular, 
$H_n\big(\Der(A, B;
{\phi})\big)$ will denote  the homology in degree $n$ 
of the chain complex
$\big(\Der_{*}(A, B; {\phi}), \delta\big)$.

The derivation complexes  enjoy  the same functorality as 
function spaces.   We will be  interested
in the chain map  $\psi^* \colon \Der_*(A, B; \phi) \to \Der_*(A', B;\phi')$
induced by pre-composition with  a DG algebra
map $\psi \colon A' \to A$ where  $\phi' = \phi \circ \psi$. 
By a standard construction, the map $\psi^*$ gives rise to a long exact homology sequence on homology 
of the form 
$$
\xymatrix{  \cdots\ar[r] & H_{n+1}\big(\Der(A', B; \phi') \big)  
\ar[r]^{\ \ J_{n+1}} & H_{n+1}\big(\Rel(\psi^*) \big)
  \ar `d[l]  `[lld]_{P_n} [lld] 
\\
H_n\big(\Der(A, B;\phi)\big) \ar[r]^{H(\psi^*)}
   & 
 H_{n}\big(\Der(A', B; \phi') \big) \ar[r] & \ldots.}
$$  
Here the ``relative'' term  $\Rel_*(\psi^*)$  is the DG space 
$\Der_*(A, B; \phi) 
\oplus  \Der_{*-1}(A', B; \phi')$ 
with differential  $D \colon \Rel_n(\psi^*) \to 
\Rel_{n-1}(\psi^*)$  defined  by the rule 
$$D(\theta, \varphi) = (\delta(\varphi) - 
\psi^*(\theta),\delta'(\theta)).$$ We have written $\delta$ and $\delta'$ for the 
differentials in $\Der_{*}(A, B; \phi)$ and $\Der_{*}(A', B; \phi'),$
respectively. The maps $J_n$ and $P_n$ are induced by the inclusion
$j_n \colon \Der_n(A, B; \phi) \to \Rel_n(\psi^*)$ and 
the projection $p_n\colon \Rel_n(\psi^*) \to \Der_{n-1}(A', B; \phi').$  

In \cite{L-S}, we showed that 
the homology theory of 
derivation complexes may be used to model
the rational homotopy theory of function
spaces at the level of the higher homotopy groups.  
The link is provided by a map 
\begin{equation} \label{eq:Phi} \Phi_f \colon \pi_n\big(\map(X,Y;f)) 
\to H_n(\Der(\mathcal{M}_Y, 
\mathcal{M}_X; \mathcal{M}_f)) \end{equation} 
whose construction originates with Thom \cite{Th}.  
Suppose $\beta \in \pi_n\big(\map(X,Y;f)\big)$
has adjoint $B \colon S^n\times X \to Y$.  Then $B$ induces a DG 
algebra map
 $\mathcal{A}_B \colon \mathcal{M}_Y \to
H^*(S^n, \Q) \otimes \mathcal{M}_X$, that is of the form
$$\mathcal{A}_B(\chi) = 1\otimes \mathcal{M}_{f}(\chi) + u\otimes
\theta_B(\chi),$$
where $u \in H^n(S^n ,\Q)$ denotes a generator  and  $\chi \in 
\mathcal{M}_X$ is of positive degree.     Here $H^*(S^n, \Q)$   
is viewed as a DG algebra with 
zero differential.   This expression
defines a linear map $\theta_B \colon \mathcal{M}_Y \to
\mathcal{M}_X$ that reduces degree by $n$.   A standard check shows 
$\theta_B$ is an $\mathcal{M}_f$-derivation cycle.  Set $\Phi_f(\beta) =
[\theta_B] \in H_n\big(\Der(\mathcal{M}_Y, \mathcal{M}_X; \mathcal{M}_f)\big)$.

\begin{proposition} \label{prop:Phi} Let $f \colon X \to Y$ be a map between
    nilpotent   CW 
complexes with $X$ finite. Then for $n \geq 2$ the map $$\Phi_f \colon \pi_n\big(\map(X,Y;f)) 
\to H_n(\Der(\mathcal{M}_Y, 
\mathcal{M}_X; \mathcal{M}_f))$$ is a well-defined natural 
homomorphism and a rational equivalence.  If $Y$ is a rational 
$H$-space,    
this holds for
$n \geq 1$ also.  Finally, in the general case for $n=1,$ $\Phi_f$ 
is a well-defined natural map of sets. 
\end{proposition}
\begin{proof} The result for $Y$ a rational $H$-space is a 
direct generalization of the identification
$ \pi_n(\map(X, K(G, n); 0) \otimes \Q  
\cong \Hom_n\big(G \otimes \Q, H^*(X, \Q)\big)$
for $G$ abelian due to Thom \cite[Th. 2]{Th}.  
The result for $n \geq 2$ is \cite[Th. 2.1]{L-S}.
The proof given there assumes $X$ and $Y$ are simply connected.
However,  all that  is needed is 
 that the map $f$ have a minimal model and that a rationalization of $Y$ 
 induces a rationalization of $\map(X, Y;f)$ as in 
 (\ref{eq:map rationalized}).
 For the case $n=1,$
 we observe that the proof that $\Phi_f$ is a well-defined  function \cite[Th. A.2]{L-S}
 for $n \geq 2$ goes through unchanged for $n =1$.  
\end{proof}
\begin{remark}
 By naturality, we mean the
following:  Suppose given maps of spaces $g\colon Y \to Z$.  Then
we have induced maps 
$$(g_*)_\sharp \colon \pi_n\big(\map(X,Y;f)\big) \to 
\pi_n\big(\map(X,Z; g\circ f)\big)$$ 
and $$H(\mathcal{M}_g^*) \colon H_n\big(\Der(\mathcal{M}_Y, \mathcal{M}_X; 
\mathcal{M}_{f}) \big) \to 
H_n\big(\Der(\mathcal{M}_Z, \mathcal{M}_X; 
\mathcal{M}_{g \circ f})\big).$$  
Then  
$H\big(\mathcal{M}_g^*\big) \circ\Phi_{f} = \Phi_{g\circ f}\circ
(g_*)_{\#}$. Note that this  identity still makes sense when $n=1$, but we must
think of both sides as maps of sets as opposed to groups.   
\end{remark}

We now return to    the geometric situation of 
\secref{sec:minimal models}. 
Let $Y$ be the total space of a principal fibration of nilpotent spaces
$K(G, n) \stackrel{j}{\to} Y \stackrel{p}{\to} B$   with   $k$-invariant $k \colon B \to 
K(G, n+1)$.  We assume $n \geq 1$ and so $G$ is abelian.  
Let $f \colon X \to Y$ be given and write $g = p \circ f \colon X \to B$. 
We have in mind a Postnikov section of $Y$ but we only need here that $k$ 
vanishes on homotopy groups. By \corref{cor:principal} we have a fibre sequence $$\xymatrix{ \map(X, Y;f) \ar[r]^{p_*} & \map(X, 
B; g) \ar[r]^{\! \! \! \! \! \! \! \! \! \! \! \! \! \! \! 
\! \! \! \! \!  k_*} & \map(X, K(G, n+1);k\circ g)}$$
and so a long exact homotopy sequence
\begin{equation} \label{eq:long homotopy}
\xymatrix@C=13pt{
&  \cdots \ar[r] &
 \pi_{n+1}\big(\map(X, K(G, n+1), k \circ g) \big)
  \ar `d[l]  `[lld]_(0.7){\partial_{n+1}} [lld] \\
\pi_{n}\big(\map(X, Y;f)\big)
\ar[rr]^{(p_*)\sharp} & &
\pi_{n}\big(\map(X, B; g)\big) \ar[r] &
\cdots 
} 
\end{equation}
$$\xymatrix{ \cdots  \ar[r]  & \pi_{1}\big(\map(X, B; g) \big) 
\ar[rr]^{\! \! \! \! \! \! \! \! \! \! \! \! \! \! \! 
\! \! \! \! \! \! \! \! \! (k_*)_\sharp} & & 
\pi_{1}\big(\map(X, K(G, n+1); k \circ g)\big)}
$$

On the level of Sullivan minimal models, we may write 
$\mathcal{M}_Y = \mathcal{M}_B \otimes \Lambda\overline{V}$ where
$\overline{V}$ is concentrated in degee $n$ and isomorphic to $G \otimes \Q.$
The differential $d_Y$ restricts to $d_B$ on the factor $\mathcal{M}_B$
while $d_Y\big(\Lambda\overline{V}\big) \subseteq \mathcal{M}_B$  is contained 
in the decomposables.  
The minimal model $\mathcal{M}_p \colon \mathcal{M}_B \to \mathcal{M}_Y$ 
may be taken to be   the inclusion.  The 
minimal model for $k$ is a map $\mathcal{M}_k \colon \Lambda V \to 
\mathcal{M}_B$ where $V$ is  isomorphic to $G \otimes \Q$ concentrated
in degree $n+1$ and $\Lambda V $ has trivial differential.
We may assume $\mathcal{M}_k(v) = d_Y(\overline{v})$
where $v \in V$  and $v \mapsto \overline{v}$ is the obvious degree
$+1$ identification of $V$ with $\overline{V}.$ 
We consider the long exact homology sequence of the map
$$\mathcal{M}_k^* \colon \Der_*\big(\mathcal{M}_B, \mathcal{M}_X; 
\mathcal{M}_g\big) \to \Der_*\big(\Lambda V, \mathcal{M}_X; 
\mathcal{M}_{k \circ g}\big).$$ 
Given a pair $(\theta, \varphi) \in \Rel_n\big( \mathcal{M}_{k}^*)$
we obtain $\overline{\theta} \in \Der_{n-1}(\mathcal{M}_Y, 
\mathcal{M}_X; \mathcal{M}_f)$ by setting 
$\overline{\theta}(\chi) =(-1)^{n-1}\varphi(\chi)$ for $\chi \in \mathcal{M}_B$
and $\overline{\theta}(\overline{v}) = \theta(v).$ 
The assignment $(\theta, \varphi) \mapsto \overline{\theta}$ 
gives a chain equivalence and thus an isomorphism
$$\Psi \colon H_n\big(\Rel(\mathcal{M}_k^*)\big) 
\stackrel{\cong}{\longrightarrow} 
 H_{n-1}\big(\Der(\mathcal{M}_Y, \mathcal{M}_X; \mathcal{M}_f)\big).$$
Observe   $\Psi \circ P_n = \mathcal{M}_{p}^*$ and write
$\Delta_n = J_n \circ \Psi$.  We  then  have a long exact sequence
\begin{equation} \label{eq:long derivation} 
\xymatrix@C=13pt{
&  \cdots \ar[r] &
 H_{n+1}\big(\Der(\Lambda V, \mathcal{M}_X; 
 \mathcal{M}_{k \circ g}) \big)
  \ar `d[l]  `[lld]_(0.7){\Delta_{n+1}} [lld] \\
H_{n}\big(\Der(\mathcal{M}_{Y}, \mathcal{M}_{X}; \mathcal{M}_f)\big)
\ar[rr]^{H(\mathcal{M}_{p}^*)} & &
H_{n}\big(\Der(\mathcal{M}_{B}, \mathcal{M}_{X};
\mathcal{M}_{g})\big) \ar[r] &
\cdots 
} 
\end{equation} 
$$\xymatrix{ \cdots  \ar[r]  & H_{1}\big(\Der(\mathcal{M}_{B},
\mathcal{M}_{X}; \mathcal{M}_g)\big) 
\ar[rr]^{H(\mathcal{M}_{k}^*)} & & 
H_{1}\big(\Der(\Lambda V, \mathcal{M}_{X}; \mathcal{M}_{k \circ g})\big)}.
$$
Applying   $\Phi$ from (\ref{eq:long homotopy}) to 
(\ref{eq:long derivation}) yields   a ladder of long exact 
sequences of groups with each vertical map a homomorphism except in the third and second 
to last terms. 
\begin{equation} \label{eq:ladder} 
    \xymatrix@C=+30pt{
 \cdots \ar[r] &    \pi_n\big(\map(X, Y; f)\big) \ar@<1ex>[d]^{\Phi_f} 
    \ar[r]^{(p_*)_\sharp}  & 
    \pi_n\big(\map(X, B; g) \big) \ar[d]^{\Phi_{g}} \ar[r] & \cdots \\
    \cdots \ar[r] &
   H_{n}\big(\Der(\mathcal{M}_{Y},
\mathcal{M}_{X}; \mathcal{M}_f)\big) 
\ar[r]^{H(\mathcal{M}_{p}^*)} &  
H_{n}\big(\Der(\mathcal{M}_{Y},
\mathcal{M}_{B}; \mathcal{M}_g)\big) \ar[r] & \cdots}     
 \end{equation}
 $$ \xymatrix@C=+30pt{
    \cdots \ar[r] & \pi_1\big(\map(X, B; g)\big) \ar[d]^{\Phi_g} 
    \ar[r]^{\! \! \! \! \!  \! \! \! \! \! \! \! \! \! \!\! \! (k_*)_\sharp}  & 
    \pi_1\big(\map(X, K(G, n+1); k \circ g) \big) \ar[d]^{\Phi_{k 
    \circ g}} \\
    \cdots \ar[r] &  H_{1}\big(\Der(\mathcal{M}_{B},
\mathcal{M}_{X}; \mathcal{M}_g)\big) 
\ar[r]^{\! \! \! \! \! \! \! \! \! H(\mathcal{M}_{k}^*)} &  
H_{1}\big(\Der(\Lambda V, \mathcal{M}_{X}; \mathcal{M}_{k \circ g})\big)}
$$
\begin{theorem}\label{thm:ladder} 
    The ladder (\ref{eq:ladder})  is commutative. 
    \end{theorem}
\begin{proof}   
We need to check commutativity of three types of squares. However, 
for two of the three, those involving $p$ and $k$,   
commutativity is a  direct consequence of the 
naturality of $\Phi$. Commutativity of the third type of square, namely 
$$ \xymatrix{
   \pi_{n+1}\big(\map(X, K(G, n+1); k \circ g)\big) \ar[d]^{\Phi_{k \circ g}} 
    \ar[r]^{\ \ \ \ \ \ \ \ \partial_{n+1}}  & 
    \pi_n\big(\map(X, Y; f) \big) \ar[d]^{\Phi_{f}} \\
   H_{n+1}\big(\Der(\Lambda V,\mathcal{M}_{X}; \mathcal{M}_{k \circ g})\big) 
\ar[r]^{\Delta_{n+1}} &  
H_{n}\big(\Der(\mathcal{M}_Y, \mathcal{M}_{X}; \mathcal{M}_{f})\big),}
$$ 
can     deduced from  naturality also, this time  using  the fibre inclusion 
$j \colon K(G, n) \to Y.$  
For,  by Thom \cite[p. 32]{Th}, in our situation  the fibre of 
the induced fibration $p_* \colon \map(X, Y;f) \to \map(X, B; g)$ may be 
identified as
the full  (possibly disconnected) function space $\map(X, K(G, n))$.   
Comparing the long 
exact homotopy sequences of the fibration 
$$\xymatrix{ \map(X, K(G, n)) \ar[r]^{j_*} & \map(X, Y;f) \ar[r]^{p_*} &
\map(X, B;g)}$$
with that of  the fibre sequence 
$$\xymatrix{ \map(X,Y) \ar[r]^{p_*} & \map(X, B;g) 
\ar[r]^{\! \! \! \! \! \! \! \!\! \! \! \!\! \! \! \!  \! \! \! k_*} &
\map(X, K(G, n+1);k \circ g)}$$
yields a factorization of $\partial_{n+1}$ of the form  
$$\xymatrix@C=-10pt{
 \pi_{n+1}\big(\map(X, K(G, n+1); k\circ g)\big) 
 \ar[rr]^{\ \ \ \ \ \ \partial_{n+1}} 
 \ar[dr]^{\cong}_{\partial'_{n+1}} & &  \pi_{n}\big(\map(X, Y; f) \big).
 \\
 &
 \pi_{n}\big( \map(X, K(G, n))\big) \ar[ur]_{(j_*)_\sharp} 
} $$
Here  $\partial'_{n+1}$ denotes the connecting homomorphism in the long exact 
homotopy sequence of the  path fibration
$\map(X, K(G, n)) \to \map(X, PK(G, n+1)) \to \map(X, K(G, n+1)).$

Such a factorization occurs, analogously,  on derivation complexes.  
Define  a  degree-lowering chain equivalence  
$  \Der_{n+1}(\Lambda V, \mathcal{M}_X; 
\mathcal{M}_{k \circ g})  
\cong  \Der_{n}(\Lambda(\overline{V}), \mathcal{M}_X; 
0)$  
  by the assignment $\theta \mapsto \overline{\theta}$  where 
 $\overline{\theta}(\overline{v}) = \theta(v)$ for $v \in V$ and 
 $\overline{\theta}$ vanishes on decomposables.
It is easy to check that,  if  
$$\Delta'_{n+1} \colon H_{n+1}\big(\Der(\Lambda V, \mathcal{M}_X; 
\mathcal{M}_{k \circ g}) \big) 
\stackrel{\cong}{\longrightarrow} H_n\big(\Der(\Lambda(\overline{V}), \mathcal{M}_X; 
0) \big)$$
is the isomorphism induced on homology, then $\partial'_{n+1} \circ 
\Phi_{k \circ g} = \Phi_{0} \circ 
\Delta'_{n+1}$. Commutativity of the third type of square now follows 
from naturality with respect to $j \colon K(G, n) \to Y.$  
\end{proof}

\section{Proof of the Main Results}\label{sec:proofs}
In this section we prove \thmref{thm:general component} 
and deduce  \thmref{thm:null component} as 
a  consequence.
\begin{proof}[Proof of \thmref{thm:general component}]
We first establish the desired formula  for an arbitrary Postnikov section
$Y^{(r)}$ of $Y$. That is, 
we prove  
\begin{equation} \label{eq:Postnikov section}
    \rank\big(\pi_1\big(\map(X, Y^{(r)}; f^r)\big)\big) = 
    \dim_\Q\big(H_1(\Der(\mathcal{M}_{Y^{(r)}}, \mathcal{M}_X; 
    \mathcal{M}_{f^r})\big) \big)
    \end{equation}
    for any given $r \geq 0$ where $f^r \colon X \to Y^{(r)}$ denotes the map
    induced by $f \colon X \to Y.$ 
    We prove this formula by induction on $r \geq 0$. 
    Since $Y^{(0)} = *,$ the base case is trivial.
    
For the induction step, we rely on  our notation 
from \secref{sec:minimal models}. 
The $r$th 
Postnikov section    $p_r\colon 
Y^{(r)} \to Y^{(r-1)}$ of $Y$ factors into a sequence of principal
fibrations
$$
Y^{(r)} = Y^{(r)}_{c_r} \to Y^{(r)}_{c_r - 1} \to \cdots \to
Y^{(r)}_{1}\to Y^{(r)}_{0},$$
with each fibration induced from a path-loop fibration 
by a map of the form $k^{r}_j
\colon  Y^{(r)}_{j-1} \to K(G^r_j, r+1)$. Let 
 $f^r_j \colon X \to Y^{(r)}_j$ denote the map induced by $f \colon X \to Y$ 
Write
$Z^{(r)}_{j}=\map(X, Y^{(r)}_{j};f^{r}_j)$. Then
\corref{cor:principal} gives a corresponding sequence of principal
fibrations
$$\map(X, Y^{(r)};f^{r}) = Z^{(1)}_{c_1} \to Z^{(1)}_{c_1 - 1}
\to \cdots \to Z^{(r)}_{1}\to Z^{(r)}_{0} = \map(X, Y^{(r-1)}; 
f^{r-1}),
$$
with each stage induced by a $k$-invariant of the form
$$(k^{r}_j)_* \colon Z^{(r)}_{j-1} \to 
\map(X, K(G^{r}_j, r+1);k^{r}_j\circ f^{r}_{j - 1}).$$
Write $ k^{(r)}_{j}= k^{r}_j\circ f^{r}_{j - 1}\colon X \to
K(G^{r}_j, r+1)$.        
The long exact homotopy sequence of the fibration
sequence
$$
Z^{(r)}_{j} \to Z^{(r)}_{j-1} \to \map(X, K(G_j, r+1);k^{(r)}_{j -1}),
$$
gives an exact sequence
$$
\xymatrix{\cdots \ar[r] & \pi_2\big(Z^{(r)}_{j-1}\big)
\ar[r]^-{k^r_{j*\#}} & \pi_2\big(\map(X,K(G_j,
r+1);k^{(r)}_{j-1})\big)
 \ar `d[l]  `[lld]_-{\partial_2} [lld]
\\
\pi_1\big(Z^{(r)}_{j}\big) \ar[r] & \pi_1\big(Z^{(r)}_{j-1}\big)
\ar[r]^-{k^r_{j*\#}} & \pi_1\big(\map(X,K(G_j,
r+1);k^{(r)}_{j-1})\big)  }
$$
This displays $\pi_1\big(Z^{(r)}_{j}\big)$ as a central extension of the 
kernel of 
$k^r_{j*\#}$ in degree one by the 
cokernel of $k^r_{j*\#}$ in degree two.
Write 
$$C^{r}_j =  \pi_2\big(\map(X,K(G_j,
r+1);k^{(r)}_{j-1})\big) / k^r_{j*\#}(\pi_2(Z^{(r)}_{j-1})) $$  
for the cokernel of $k^r_{j*\#}$ in degree $2$ and 
$$I^{r}_{j} = \im\{ k^r_{j*\#}\colon \pi_1(Z^{(r)}_{j-1}))
\to \pi_1\big(\map(X,K(G_j,r+1);k^{(r)}_{j-1})\big) \}$$ for its image  on fundamental groups. Then  
\begin{equation}
    \label{eq:rank homotopy}
    \rank(\pi_1\big(Z^{(r)}_j\big)\big)   = \,  
\rank(\pi_1\big(Z^{(r)}_{j-1}\big)\big)
  \,  + \, \rank(C^{r}_{j}) \, - \, \rank(I^{r}_j).
\end{equation}

We can follow  the preceding line of reasoning, analogously,  within  the framework of 
derivation complexes of minimal models.  Write
$D_j^{(r)} = \Der_*(\mathcal{M}_{Y^{(r)}_j}, \mathcal{M}_X; 
\mathcal{M}_{f^{r}_{j}}).$ We have a sequence
$$\begin{array}{lll} 
\Der_*(\mathcal{M}_{Y^{(r)}}, \mathcal{M}_X; \mathcal{M}_{f^r}) & =
D^{(r)}_{c_r} \to    D^{(r)}_{c_r - 1} \to \cdots & \\ 
& \\  & \cdots \to
D^{(r)}_1 \to D^{(r)}_0 = 
\Der_*(\mathcal{M}_{Y^{(r-1)}}, \mathcal{M}_X; \mathcal{M}_{f^{r-1}}).
& \end{array}
$$

The long exact homology sequence corresponding to  the map 
$$ \mathcal{M}_{k^{r}_j} \colon D^{(r)}_j \to \Der_*(\Lambda V_{j}^{(r)}, \mathcal{M}_X; 
\mathcal{M}_{k^{(r)}_j}) $$
takes the form
$$
\xymatrix@C=+38pt{\cdots \ar[r] & H_2\big(D^{(r)}_{j-1}\big)
\ar[r]^-{\! \! \! \!  H(\mathcal{M}_{k^r_{j}}^*)} &  
H_2\big(\Der(\Lambda V_{j}^{(r)}, \mathcal{M}_X; 
\mathcal{M}_{k^{(r)}_j})\big)
 \ar `d[l]  `[lld]_-{\Delta_2} [lld]
\\
H_1\big(D^{(r)}_{j}\big) \ar[r] & H_1\big(D^{(r)}_{j-1}\big)
\ar[r]^-{\! \! \! \! H(\mathcal{M}_{k^r_{j}}^*)} & H_1\big(\Der(\Lambda V_{j}^{(r)}, \mathcal{M}_X; 
\mathcal{M}_{k^{(r)}_j})\big)  }
$$
Writing 
$$\mathcal{C}^{r}_j =  H_2\big(\Der(\Lambda V_{j}^{(r)}, \mathcal{M}_X; 
\mathcal{M}_{k^{(r)}_j})\big) / H(\mathcal{M}_{k^r_{j}}^*)(H_2(D^{(r)}_{j-1})) $$  
for the cokernel and 
$$\mathcal{I}^{r}_{j} = \im\{ H(\mathcal{M}_{k^r_{j}}^*)\colon H_1(D^{(r)}_{j-1})
\to H_1\big(\Der(\Lambda V_{j}^{(r)}, \mathcal{M}_X; 
\mathcal{M}_{k^{(r)}_j})\big) \}$$ for  the image 
of
$H(\mathcal{M}_{k^r_{j}}^*)$ as above, we obtain   
\begin{equation}
    \label{eq:rank derivation}
    \dim_\Q(H_1\big(Z^{(r)}_j\big)\big)   = \,  
\dim_\Q(H_1\big(Z^{(r)}_{j-1}\big)\big)
  \,  + \, \dim_\Q(\mathcal{C}^{r}_j) \, - \, 
  \dim_\Q(\mathcal{I}^{r}_j).
\end{equation}

Now, combining \propref{prop:Phi} and \thmref{thm:ladder} we
see $$\rank(C_j^r) = \dim_\Q(\mathcal{C}^r_j) \hbox{\ \ and \ \ }
\rank(I_j^r) = \dim_\Q(\mathcal{I}^r_j).$$
An easy induction on 
$j = 1, \ldots, c_r$ using (\ref{eq:rank homotopy}) and (\ref{eq:rank derivation})
completes the main induction step and establishes 
(\ref{eq:Postnikov section}) for all $r \geq 0.$  

To complete the proof of \thmref{thm:general component}
we observe that, by obstruction theory,  the $r$-equivalence $p_r \colon Y \to
Y^{(r)}$ induces an $(r-N)$-equivalence $(p_r)_* \colon
\map(X,Y;f) \to \map(X,Y^{(r)};f^{r})$ where $N$ is the dimension 
of the finite complex $X$. 
On the derivation complex side,  it is straightforward to prove  
$$\mathcal{M}_{p_r}^{*} \colon \Der_*(\mathcal{M}_{Y^{(r)}}, 
\mathcal{M}_X; \mathcal{M}_{f^{r}}) \to
\Der_*(\mathcal{M}_{Y}, 
\mathcal{M}_X; \mathcal{M}_f)
$$
induces an $(r-N)$-homology equivalence.  Thus, for $ r \geq N+1,$ 
$$\begin{array}{ll} \rank(\pi_1\big(\map(X, Y; f)) & =
\rank\big(\pi_1\big(\map(X, Y^{(r)}; f^{r})\big)\big) \\
& = \dim_\Q\big(H_1\big(\Der(\mathcal{M}_{Y^{(r)}}, \mathcal{M}_X; 
\mathcal{M}_{f^{r}})\big)\big) \\
& = \dim_\Q\big(H_1\big(\Der(\mathcal{M}_{Y}, \mathcal{M}_X; 
\mathcal{M}_f)\big)\big)
\end{array}
$$
\end{proof}
\begin{proof}[Proof of \thmref{thm:null component}]

Observe that 
$$H_1\big(\Der(\mathcal{M}_Y, \mathcal{M}_X; 0)\big)
\cong \Hom_1\big(Q(\mathcal{M}_Y), H(\mathcal{M}_X)\big)$$
where we recall $Q(\mathcal{M}_Y)$ denotes the quotient space 
of indecomposables of the minimal model of $Y$.  The result
now follows directly from the fact that  
$Q_n(\mathcal{M}_Y) \cong \pi_n(Y)_\Q$ for $n \geq 2$
while $\dim_\Q(Q_1(\mathcal{M}_Y)) = \rank(\pi_1(Y))$
by the results mentioned in \secref{sec:minimal models}.
\end{proof}
\section{Consequences and  Examples}\label{sec:consequences}
First we  illustrate   the role of the 
map $f \colon X \to Y$ for the rank of the fundamental group 
of $\map(X, Y;f).$  
As a direct consequence of \thmref{thm:general component}, we obtain 
the following basic fact.  
\begin{theorem} \label{thm:inequality} 
    Let $f \colon X \to Y$ be a map between nilpotent CW complexes
    of finite type
    with $X$ finite.  Then   
    $$ \rank\big(\pi_1\big(\map(X, Y; f)\big)\big) \leq 
    \rank\big(\pi_1\big(\map(X, Y;0)\big)\big).$$
\end{theorem}
\begin{proof}   It suffices to note
    $\dim_\Q\big(H_n\big(\Der(A, B; \phi) \big) \big)  
    \leq \dim_\Q\big(H_n\big(\Der(A, B; 0) \big) \big)$
    for any DG algebra map $\phi \colon A \to B.$ 
\end{proof}

The fundamental group can  distinguish components   of $\map(X, Y)$   
even in a simple case such as    $Y= K(G, 1)$ for $G$ non-abelian. For in this case,  
 by a result of Gottlieb,  $\map(X, K(G, 1)) \simeq K(C(f_\sharp), 1)$ 
where $C(f_\sharp)$ is the centralizer of the map $f_\sharp \colon 
\pi_1(X) \to G$ 
\cite[Lem.2]{Got3}. 
In particular, Gottlieb's result describes the structure of the fundamental group  of $\map(X, Y;f)$ 
directly in terms of     
$f_\sharp$. 
We complement this result  with the following
example.   
\begin{example} \label{ex:free loop}  Let $Y$ be a nilpotent CW complex
   and  $f \colon S^1 \to Y$.   We compute the rank of the fundamental group of 
   the component $\map(S^1, Y;f)$ of the free loop space in terms of the 
   rational homotopy class
     $\alpha \in \pi_1(Y)_\Q$ represented by $f$. 
   If $\pi_2(Y) \otimes \Q = 0$ then, by (\ref{eq:map rationalized}) 
   and the argument of Hansen 
   \cite[Prop.1]{Han}, $\pi_1\big(\map(S^1, Y; f)\big)_\Q = 
   C(\alpha)$  
   where $C(\alpha)$ denotes 
    the centralizer of $\alpha$ in $\pi_1(Y)_\Q.$ In general, we prove  
    \begin{equation}\label{eq:centralizer}  \rank\big(\pi_1\big(\map(S^1, Y; f)\big)\big) 
    = \rho_2(Y) + \rank(C(\alpha)). \end{equation}
  
    To begin, as in the last step of the proof 
  of \thmref{thm:general component}, we have 
   $$H_1\big(\Der(\mathcal{M}_Y, \mathcal{M}_{S^1}; \mathcal{M}_f)\big)
  \cong  H_1\big(\Der(\mathcal{M}_{Y^{(2)}}, 
  \mathcal{M}_{S^1}; \mathcal{M}_{f^2})\big)
 $$ where $\mathcal{M}_{Y^{(2)}}$ is the $2$-minimal model for $Y$
 and $f^2 \colon S^1 \to Y^{(2)}$ the induced map. 
 We compute the dimension of the latter space.    Fixing  notation as  
  in \secref{sec:minimal models},   write 
 $(\mathcal{M}_{Y^{(2)}},d_{Y^{(2)}}) = 
 (\Lambda \big( W^{(1)} \oplus W^{(2)} \big), D)
 $
 where $W^{(m)} = \overline{V}^{(m)}_1 \oplus \cdots \oplus 
 \overline{V}^{(m)}_{c_m}$ for $ 
 m=1,2.$ Recall   $D(\overline{V}^{(m)}_1) = 0$ while  
  $D(\overline{V}^{(m)}_j) \subseteq \Lambda (W^{(m-1)} \oplus 
 \overline{V}^{(m)}_1 \oplus \cdots \oplus \overline{V}^{(m)}_{j-1})$ for $j =2, 
 \ldots, c_m,$ $m= 1,2$  and where $W^{(0)} = \{0\}$.
 Fix bases $\overline{V}_{j}^{(m)} = \Q(\overline{v}^{(m)}_{j, 1}, 
 \ldots, \overline{v}^{(m)}_{j,n^m_j}).$ 
 Write $\mathcal{M}_{S^1} = \Lambda(t)$ for $t$ of degree one.
 Without loss of generality, 
 we assume     the 
homotopy class $\alpha$ corresponds to a basis element 
$\overline{v}^{(1)}_{j_0, k_0}$ via Sullivan's isomorphism.
That is, we assume $\mathcal{M}_{f^2}(\overline{v}^{(1)}_{j_0, k_0}) =t$
while $\mathcal{M}_{f^2}(\overline{v}^{(1)}_{j, k}) =0$
for $(j, k) \neq (j_0, k_0).$ 
The centralizer $C(\alpha)$, in this  set-up,
corresponds to the space 
$C(\overline{v}^{(1)}_{j_0, k_0})$  spanned by those basis vectors 
$\overline{v}^{(1)}_{j, k}$ such  
that no non-zero multiple of the product 
$\overline{v}^{(1)}_{j_0, k_0} \cdot \overline{v}^{(1)}_{j, k}$  
appears as a  summand in $D(\overline{v}^{(1)}_{j', k'})$ for any $j', 
k'.$ 

We observe that there are no boundaries in degree one in 
 the chain complex $\Der_*(\mathcal{M}_{Y^{(2)}}, 
  \mathcal{M}_{S^1}; \mathcal{M}_{f^2}).$
  For given an $\mathcal{M}_{f^2}$-derivation $\theta$ of  degree $2$
 and $\chi \in \mathcal{M}_{Y^{(2)}},$ we have 
 $\delta(\theta)(\chi) =  (-1)^n\theta(D(\chi)) 
 =0$ since $d_{S^1}$ is trivial while $D(\chi)$ is decomposable
 and $\mathcal{M}_{f^2}$ vanishes above degree $1$.
 
As for degree one cycles, write $\theta_{j,k}$ for
the $\mathcal{M}_{f^2}$-derivation which carries
 $\overline{v}^{1}_{j,k}$ to $1 \in \mathcal{M}_{S^1}$ and vanishes on
 the other basis elements of $\mathcal{M}_{Y^{(2)}}.$ Similarly, let
 $\varphi_{j,k}$ carry $\overline{v}^2_{j,k}$ to $t$ and all other
basis elements to $0.$  It is easy to check that $\delta 
(\varphi_{j,k}) = 0$ for all $j,k$ which accounts for the $\rho_2(Y)$
term above. We also see that if $v_{j,k}^{(1)} \in 
Z(\overline{v}^{(1)}_{j_0, k_0})$
then $\delta(\theta_{j,k}) =0$. Finally, suppose $\overline{v}_{j, 
k}^{(1)} \not\in Z(\overline{v}^{(1)}_{j_0, k_0})$ for some $j,k.$
Then we may choose $\overline{v}^{(1)}_{j', k'}$ such that 
$$D\big(\overline{v}^{(1)}_{j', k'}\big) = 
c\cdot \overline{v}^{(1)}_{j_0, k_0} \cdot \overline{v}^{(1)}_{j, k} \ \
+\hbox{ \ \ other quadratic terms }$$ 
for $c \neq 0$.  We then see $$\delta\big(\theta_{j, k}\big)
\big(\overline{v}^{(1)}_{j', k'}\big) = 
\theta_{j, k}\big(D(\overline{v}^{(1)}_{j', k'})\big) =
c \cdot \mathcal{M}_{f^2}(v^{(1)}_{j_0, k_0}) \cdot 
\theta_{j, k}\big(\overline{v}^{(1)}_{j, k}\big) 
= c \cdot t \neq 0,$$ and so $\theta_{j, k}$
is not a cycle. This establishes (\ref{eq:centralizer}).  
 \end{example}

 We next consider a class of examples 
 for which the   rank of the fundamental
 group of $\map(X, Y;f)$   depends only on the  map 
 $H(f) \colon H^*(Y, \Q) \to H^*(X, \Q).$   
  We say   a simply connected CW complex $Y$ is an {\em $F_0$-space}
    if $Y$ is rationally  
    elliptic (rational homology and rational homotopy both 
    finite-dimensional) with  
    positive Euler characteristic.  Equivalently, an $F_0$-space is 
    any elliptic complex with vanishing  rational cohomology in odd degrees.  
    Examples of $F_0$-spaces include 
    (products of) even dimensional spheres, complex projective spaces 
    and, more generally,  
    any homogeneous spaces $G/H$ with $H$ a closed subgroup of maximal rank.   Following 
    Grivel \cite{Griv}, we can compute the rank of  
    $\pi_1\big(\map(X, Y;f)\big)$ for $f \colon X \to Y$ a map between 
    $F_0$-spaces   directly in terms of  the degree $2$    
    cohomology derivation  space  $\Der_2(H^*(Y, \Q), H^*(X, \Q); H(f))$.  
    \begin{theorem}\label{thm:F_0-space}
    Let $f \colon X \to Y$ be  a map 
    between $F_0$-spaces where $H^*(Y, \Q)$ has top degree $2N$. Let     
    $D_{2}(f) = \dim_\Q\big( \Der_{2}(H^*(Y, \Q), H^*(X, \Q); 
    H(f))\big)$. 
    Then
    $$  \rank(\pi_1(\map(X,Y;f)) = D_2(f) + \sum_{i=1}^{N}\rho_{2i+1}(Y)\cdot 
b_{2i}(X) - \sum_{i=0}^{N} \rho_{2i+2}(Y)b_{2i}(X)
.$$
\end{theorem}    
  \begin{proof}
    By  results of Halperin \cite{Hal}, the minimal model
     for  $Y$ takes the form 
    $(\mathcal{M}_Y, d_Y)$ $= (\Lambda V_0\otimes
    \Lambda V_1, d_Y)$ where $V_0$ and $V_1$ are graded spaces of 
    equal (finite) dimension with $V_0$  evenly graded,
      $V_1$   oddly graded.  
      The differential satisfies  $d_Y(V_0) = \{0\}$
      while $d_Y$ maps $V_1$ into the decomposables of $\Lambda 
      V_0.$  Write
 $\rho_Y \colon \Lambda V_0 \otimes \Lambda V_1 \to H^*(Y, \Q)$ for the map
   which sends elements of 
  $V_0$ to their corresponding cohomology class and elements of $V_1$ 
  to zero. The map  $\rho_Y$  then represents a 
  formalization
  of $(\mathcal{M}_Y, d_Y).$

  The needed result  follows directly from  \thmref{thm:general component} 
  and the existence of
 an  exact sequence of the form 
   $$ 
\xymatrix{0 \to \Der_2(H^*(Y, \Q),H^*(X,\Q); 
H(f)) \ar[r]^{\rho_Y^*} &   
\Der_2(\Lambda V_0, H^*(X, \Q); 
H(f) \circ \rho_Y)
 \ar `d[l]  `[ld]_-{\! \! \! \! \! \! \! \! \! \! \! \! d_Y^*} [ld]
\\
\Der_1(\Lambda V_1, H^*(X,\Q);0) \ar[r]^{\mathcal{H}} & 
H_1\big(\Der(\mathcal{M}_Y, \mathcal{M}_X; \mathcal{M}_f)\big)
\to 0. &  }
$$
The latter  represents an extension of Grivel's   \cite[Th.4.4]{Griv} 
to our framework of generalized derivations. 
 Here $\rho_Y^*$ and $d_Y^*$ are the  maps induced by pre-composition with   
 the formalization   and the differential of the minimal model of $Y$. 
  Given   
    $\theta \in \Der_1(\Lambda V_1, H^*(X, \Q);0)$
    define 
    $H(\theta)   \in \Der_1(\mathcal{M}_Y, \mathcal{M}_X; \mathcal{M}_f),$ 
    as follows.  
   Set  $H(\theta)(x) = 0$ for $x \in V_0 $ and   
    $H(\theta)(y) = P$  for $y \in V_1$
      where $P \in \mathcal{M}_X$ is any cycle representative of the class 
     $\theta(y) \in \mathcal{M}_X, y \in V_1$.  Extend   by the 
     $\mathcal{M}_f$-derivation law.  The result is a cycle derivation
     $H(\theta)$. 
     The map $\mathcal{H} \colon \Der_1(\Lambda V_1, H^*(X, \Q);0) \to 
    H_1\big(\Der(\mathcal{M}_Y, \mathcal{M}_X; \mathcal{M}_f)\big)$
    is then   defined to be the linear map carrying $\theta$ to 
        the homology class of 
     $H(\theta).$ 
Grivel's proof  of the above cited result is   directly adapted to show 
$\mathcal{H}$ is well-defined and 
the sequence is exact. 
\end{proof}

 As a corollary, we deduce, for instance, the following example.     
    \begin{example}  \label{ex:G/T}
    Let $G$ be a compact simple Lie group of rank $n > 1$ and   
    $T^n \subseteq G$    a maximal torus. Let $f \colon S^2 \to G/T^n$ 
    be any given 
    map. 
   The space $G/T^n$ is an $F_0$-space with rational cohomology generated in 
   degree $2$. It is classical that $\rho_3(G/T^n) =1.$ Thus 
   by \thmref{thm:F_0-space}, we have 
   $$\rank\big(\pi_1(\map(S^2, G/T^n; f))\big) = 1 - n + 
   D_2(f).$$ 
   If $f$ is rationally trivial, then $D_2(f) =n$
   and so $\rank\big(\pi_1(\map(S^2, G/T^n; f))\big) = 1$. 
   
 Now suppose $f$ is rationally essential.  Fix an additive basis $\{ t_1, 
 \ldots, t_n \}$ for   $H^2(Y, \Q)$   and suppose,  say, 
   $H(f)(t_i) \neq 0$.  Let  $\theta_i$ be dual to $t_i$ in 
   the space $\Hom_2(H^*(Y, \Q), H^*(X, \Q))$: that is,     $\theta_i(t_i) = 1$ 
   while $\theta_i(t_j) = 0$ for $j \neq i.$  Suppose $\theta_i$ 
   extends to an  
   $H(f)$-derivation.  Since the Weyl group of $G$ is a finite 
   reflection group, the cohomology class  
   $t_1^2 + \cdots +
   t_n^2 $ vanishes in $H^4(G/T^n, \Q)$. On the other hand,
   $$\theta_i(t_1^2 + \cdots + t_n^2) = 2 H(f)(t_i) \theta_i(t_i)  = 2 H(f)(t_i)\neq 
   0.$$  This contradiction implies $D_2(f) = n - 1$ and 
   so $\pi_1(\map(S^2, G/T^n;f))$ is a finite group in this case. 
   
   Finally, by  the  
   result of E. Cartan, $\pi_2(G) = 0$ and so  
   $\pi_2(G/T^n)$ is free abelian.  Thus    $f \colon S^2 \to 
   G/T^n$ is essential if and only if it is  rationally essential.
   Summarizing, we have shown 
    \begin{equation} \label{eq:finite}
	\pi_1\big(\map(S^2, G/T^n; f)\big)  \hbox{  is a finite group if and 
    only if   }  f 
   \hbox{  is    essential.} \end{equation}
   
  \end{example}   
  
We conclude with some results concerning further structure
of the rationalized fundamental group of function space components. 
First we observe that, in the case of the null component, 
our inductive procedure for computing  rank  can also be applied to 
analyze the nilpotency class, and the ranks of successive
quotients in the lower central series.
For instance, we can phrase the
following result.

\begin{theorem}
Let $X$ and $Y$ be  nilpotent CW  complexes with $X$ finite 
of dimension $N$.   Suppose the Postnikov decomposition of 
$Y$ admits a principal refinement of length $c_r$ at the $r$th stage.
Then, for each $r \geq 1$, the $r$th stage of the Postnikov
decomposition of $\map(X,Y;0)$ admits a principal 
refinement of length $\leq \sum_{j=r}^{N+r} c_j$.
\end{theorem}
\begin{proof} Each principal fibration $K(G^{r}_{j}, r) \to 
Y^{(r)}_j \to Y^{(r)}$ in a Postnikov decomposition of $Y$ 
leads to a principal fibration
$\map(X, Y^{(r)}_j, 0) \to \map(X, Y^{(r)}_{j-1}; 0)$ by 
\corref{cor:principal}.  Since we are working with 
null-components, the fibre can be identified with
$\map(X, K_1; 0).$ The result now follows   
by   induction,  as  in the proof of \thmref{thm:general component}. \end{proof}

Recall a space $Y$ is {\em simple} if $\pi_1(Y)$ is 
abelian and acts trivially on the higher homotopy groups of $Y$.
In this case, $c_r=1$ for all $r$. Thus we have 
\begin{corollary}
Suppose $X$ is a nilpotent finite complex and $Y$ is a \emph{simple} space. 
Then the Postnikov decomposition of $\map(X,Y;0)_\Q$ admits a principal
refinement that is of length $\leq N$ at each stage.
In particular, $\pi_1\big(\map(X,Y;0)\big)_\Q$
is of nilpotency class $\leq N$. 
\end{corollary}

 Finally, we have seen above that the space of maps into  a (simply 
 connected) two-stage 
 Postnikov piece
$Y$ can have non-abelian rationalized fundamental group.
We conclude with a simple result going the other way. Say a nilpotent space 
$Y$ is a {\em rational two-stage} space if  
its rationalization $Y_\Q$ appears as  the total space
    of a principal fibration $K_1 \to Y_\Q \stackrel{p}{\to} K_0$ with
    $K_0$ and $K_1$ rationalized $H$-spaces.  Many spaces of interest 
    are rational two-stage including, for instance,  homogeneous 
    spaces $G/H$ for $G$ and $H$ compact Lie groups with $H$ now of 
    arbitrary rank. Write 
    $H^*(K_i,\Q) = \Lambda W_i$   for graded
    rational vector spaces $W_i$,  $i = 0,1$. 
    We have 
\begin{theorem} Let $X$ be a finite, nilpotent CW complex and $Y$ a 
rational two-stage space as above.    Suppose  $\Hom_1(W_0, H^*(X, \Q))$ and 
    $\Hom(W_1, H^*(X,\Q))$ are both trivial. Then the rationalization of 
    $\map(X, Y;f)$ is a simple space and, in particular, the group
    $\pi_1\big(\map(X, Y;f)\big)_\Q$ is abelian for all $f \colon X \to Y$.
    \end{theorem}
    \begin{proof} In this case, we obtain 
	$p_* \colon \map(X, Y_\Q; f_\Q) \to \map(X, K_0;  p \circ f_\Q),$
	a principal fibration with fibre homeomorphic to the full function 
	space $\map(X, K_1)$ as in the proof of \thmref{thm:ladder}
	above.  However,  
	$[X, K_1] \cong \Hom(W_1, H^*(X, \Q)) = 0$ by hypothesis
	and so $\map(X, K_1)$ 
	is actually connected.  Thus $p_*$ is a principal  fibration
	expressing $\map(X, Y_\Q; f_\Q)$ as a generalized two-stage
	rational space and determining a (possibly non-minimal) 
	model for $Y$.
	Again, by hypothesis, 
	$\pi_1\big(\map(X, K_0; p \circ f_\Q)\big) \cong 
	\Hom_1(W_0,H^*(X, \Q)) = 0.$ This means the differential
	in the  induced model for $\map(X, Y;f)$ has  linear
	or trivial differential in degree one. It follows that
	the minimal model for 
	$\map(X, Y;f)$ has trivial differential in degree one. 
	\end{proof}
	
We deduce directly the following result which, combined with 
\thmref{thm:F_0-space} above, gives the full
structure the rationalized fundamental group of $\map(X, Y;f)$ for 
$X$ and $Y$ $F_0$-spaces. 
  \begin{corollary} \label{cor:abelian}  Let $X$ and $Y$ be $F_0$-spaces. Then  
 all components of $\map(X, Y)$ 
 have abelian  rationalized fundamental group.  
 \end{corollary}
\providecommand{\bysame}{\leavevmode\hbox to3em{\hrulefill}\thinspace}
\providecommand{\MR}{\relax\ifhmode\unskip\space\fi MR }
\providecommand{\MRhref}[2]{%
  \href{http://www.ams.org/mathscinet-getitem?mr=#1}{#2}}
  \providecommand{\href}[2]{#2}

\end{document}